\documentclass[12pt, final]{amsart}

\usepackage{mathrsfs}
\usepackage{hyperref}
\usepackage[svgnames]{xcolor}
\hypersetup{colorlinks,breaklinks,
            linkcolor=DarkBlue,urlcolor=DarkBlue,
            anchorcolor=DarkBlue,citecolor=DarkBlue}
\usepackage{url}
\usepackage{euscript}
\usepackage[verbose,letterpaper,tmargin=1in,bmargin=1in,lmargin=1in,rmargin=1in]{geometry}
\usepackage{amssymb}

\newtheorem{theorem}{Theorem}[section]
\newtheorem{maintheorem}{Theorem}

\newtheorem{proposition}[theorem]{Proposition}
\newtheorem{corollary}[theorem]{Corollary}
\newtheorem{lemma}[theorem]{Lemma}

\theoremstyle{definition}
\newtheorem{definition}[theorem]{Definition}
\newtheorem*{acknowledgement}{Acknowledgement}

\newtheorem{remark}[theorem]{Remark}
\newtheorem{remark*}[theorem]{}

\newcommand{\C}{\mathbb C}
\newcommand{\Z}{\mathbb Z}
\newcommand{\N}{\mathbb N}

\newcommand{\mc}[1]{\mathcal{#1}}
\newcommand{\mb}[1]{\mathbb{#1}}

\renewcommand{\labelenumi}{(\arabic{enumi})}

\begin{document}

% Set the beginning of a LaTeX document
\title{A characterization of freeness by invariance under quantum spreading}
\author{Stephen Curran}
\address{Department of Mathematics, University of California, Los Angeles, CA 90095, USA.}
\email{\href{mailto:curransr@math.ucla.edu}{curransr@math.ucla.edu}}
\urladdr{\href{http://www.math.ucla.edu/~curransr}{http://www.math.ucla.edu/~curransr}}

\begin{abstract}
We construct spaces of quantum increasing sequences, which give quantum families of maps in the sense of So{\l}tan.  We then introduce a notion of quantum spreadability for a sequence of noncommutative random variables, by requiring their joint distribution to be invariant under taking quantum subsequences.  Our main result is a free analogue of a theorem of Ryll-Nardzewski: for an infinite sequence of noncommutative random variables, quantum spreadability is equivalent to free independence and identical distribution with respect to a conditional expectation.
\end{abstract}

\maketitle

\section*{Introduction}

The study of random objects with distributional symmetries is an important subject in modern probability.  Consider a sequence $(\xi_1,\xi_2,\dotsc)$ of random variables.  Such a sequence is called \textit{exchangeable} if its distribution is invariant under finite permutations, and \textit{spreadable} if it is invariant under taking subsequences, i.e., if
\begin{equation*}
 (\xi_1,\dotsc,\xi_k) \stackrel{d}{\sim} (\xi_{l_1},\dotsc,\xi_{l_k})
\end{equation*}
for all $k \in \N$ and $l_1 < \dotsb < l_k$.  In the 1930's, de Finetti gave his famous characterization of infinite exchangeable sequences of random variables taking values in $\{0,1\}$ as conditionally i.i.d.  This was extended to variables taking values in a compact Hausdorff space by Hewitt and Savage \cite{hs}.  It was later discovered by Ryll-Nardzewski that de Finetti's theorem in fact holds under the apparently weaker condition of spreadability \cite{rn}.  For a comprehensive treatment of distributional symmetries in classical probability, the reader is referred to the recent text of Kallenberg \cite{kal}.

Free probability, developed by Voiculescu in the 1980's, is based on the notion of \textit{free independence} for random variables with the highest degree of noncommutativity.  Remarkably, there is a deep parallel between the theories of classical and free probability.  However, it is only quite recently that this parallel has been extended to the study of distributional symmetries.  The breakthrough came with the work of K\"{o}stler and Speicher \cite{ksp}, who discovered that, roughly speaking, in free probability one should consider \textit{quantum} distributional symmetries.  More specifically, they defined the notion of \textit{quantum exchangeability} for a sequence $(x_1,x_2,\dotsc)$ of noncommutative random variables by requiring that for each $n \in \N$, the joint distribution of $(x_1,\dotsc,x_n)$ is invariant under the natural action of the \textit{quantum permutation group} $A_s(n)$ of Wang \cite{wang2}.  They then gave a free analogue of de Finetti's theorem: for an infinite sequence of noncommutative random variables, quantum exchangeability is equivalent to free independence and identical distribution with respect to a conditional expectation.  This has since been extended to more general sequences \cite{cur3}, and to sequences invariant under actions of other compact quantum groups \cite{cur4,bcs2}.  (See also \cite{kos} for a detailed analysis of exchangeability and spreadability for sequences of noncommutative random variables).

The purpose of the present paper is to develop a notion of \textit{quantum spreadability} for sequences of noncommutative random variables.  The first problem is to find a suitable quantum analogue of an increasing sequence.  The answer which we suggest here is similar to Wang's notion of a quantum permutation.  For natural numbers $k \leq n$ we construct certain universal C$^*$-algebras $A_i(k,n)$, which we call \textit{quantum increasing sequence spaces}, whose spectrum is naturally identified with the space of increasing sequences $1 \leq l_1 < \dotsb < l_k \leq n$.  These objects form \textit{quantum families of maps}, in the sense of So{\l}tan \cite{soltan}, from $\{1,\dotsc,k\}$ into $\{1,\dotsc,n\}$.  Quantum spreadability is naturally defined as invariance under these familes of quantum transformations.  This approach is justified by our main result, which is a free analogue of the Ryll-Nardzewski theorem for quantum spreadable sequences (see Sections \ref{background} and \ref{qinvariant} for definitions and motivating examples):

\begin{maintheorem}\label{mainthm}
Let $(\rho_i)_{i \in \N}$ be an infinite sequence of unital $*$-homomorphisms from a unital $*$-algebra $C$ into a tracial W$^*$-probability space $(M,\tau)$.  Then the following are equivalent:
\begin{enumerate}
\renewcommand{\labelenumi}{(\roman{enumi})} 
\item $(\rho_i)_{i \in \N}$ is quantum exchangeable.
\item $(\rho_i)_{i \in \N}$ is quantum spreadable.
\item $(\rho_i)_{i \in \N}$ is freely independent and identically distributed with respect to the conditional expectation $E$ onto the tail algebra
\begin{equation*}
 B = \bigcap_{n \geq 1} W^*\bigl(\{\rho_i(c): c \in C, i \geq n\}\bigr).
\end{equation*}

\end{enumerate}

\end{maintheorem}

The equivalence of (i) and (iii) is the main result of \cite{ksp} in the case $C = \C[t]$, and was shown for general $C$ in \cite{cur3}.  

Observe that Theorem \ref{mainthm} holds only for infinite sequences.  In \cite{cur3}, we have given an approximation to how far a finite quantum exchangeable sequence is from being free with amalgamation.  As in the classical case, finite quantum spreadable sequences are more difficult, and we will not attempt an analysis here.  For a treatment of classical finite spreadable sequences, see \cite{kal2}.

Our paper is organized as follows.  Section \ref{background} contains notations and preliminaries.  We recall the basic notions from free probability, and introduce Wang's quantum permutation group $A_s(n)$.  In Section \ref{qinc}, we introduce the algebras $A_i(k,n)$ and prove some basic results.  In particular we show that $A_i(k,n)$ is a quotient of $A_s(n)$.  In Section \ref{qinvariant}, we introduce the notions of quantum exchangeability and spreadability, and prove the implications (i) $\Rightarrow$ (ii) and (iii) $\Rightarrow$ (i) of Theorem \ref{mainthm}.  These implications hold in fact for finite sequences, and in a purely algebraic context.  We complete the proof of Theorem \ref{mainthm} in Section 4, by showing the implication (ii) $\Rightarrow$ (iii).

\section{Background and notations}\label{background}

\begin{remark*} \textbf{Notations}.  
Let $C$ be a unital $*$-algebra.  Given an index set $I$, we let 
\begin{equation*}
 C_I = \mathop{\ast}_{i \in I} C^{(i)}
\end{equation*}
denote the free product (with amalgamation over $\C$), where for each $i \in I$, $C^{(i)}$ is an isomorphic copy of $C$.  For $c \in C$ and $i \in I$ we denote the image of $c$ in $C^{(i)}$ as $c^{(i)}$. The universal property of the free product is that given a unital $*$-algebra $A$ and a family $(\rho_i)_{i \in I}$ of unital $*$-homomorphisms from $C$ to $A$, there is a unique unital $*$-homomorphism from $C_I$ to $A$, which we denote by $\rho$, such that $\rho(c^{(i)}) = \rho_i(c)$ for $c \in C$ and $i \in I$.  We will mostly be interested in the case that $I = \{1,\dotsc,n\}$, in which case we denote $C_I$ by $C_n$, and $I = \N$ in which case we denote $C_I = C_\infty$.
\end{remark*}

\begin{remark*}\textbf{Free Probability.} We begin by recalling some basic notions from free probability, the reader is referred to \cite{vdn},\cite{ns} for further information.
\end{remark*}

\begin{definition}\hfill
 \begin{enumerate}
  \item A \textit{noncommutative probability space} is a pair $(A,\varphi)$, where $A$ is a unital $*$-algebra and $\varphi$ is a state on $A$.
\item A \textit{W$^*$-probability space} is a pair $(M,\tau)$, where $M$ is a von Neumann algebra and $\tau$ is a faithful normal state which is tracial, i.e., $\tau(xy) = \tau(yx)$ for $x,y \in M$.
 \end{enumerate}

\end{definition}

\begin{definition}
Let $C$ be a unital $*$-algebra, $(A,\varphi)$ a noncommutative probability space and $(\rho_i)_{i \in I}$ a family of unital $*$-homomorphisms from $C$ to $A$.  The \textit{joint distribution} of the family $(\rho_i)_{i \in I}$ is the state $\varphi_\rho$ on $C_I$ defined by $\varphi_\rho = \varphi \circ \rho$.  $\varphi_\rho$ is determined by the \textit{moments}
\begin{equation*}
\varphi_\rho(c_1^{(i_1)}\dotsb c_k^{(i_k)}) = \varphi(\rho_{i_1}(c_1)\dotsb \rho_{i_k}(c_k)),
\end{equation*}
where $c_1,\dotsc,c_k \in C$ and $i_1,\dotsc,i_k \in I$.
\end{definition}

\begin{remark*}\label{motivation}\textit{Examples.}
\begin{enumerate}
\item   Let $(\Omega,\mc F, P)$ be a probability space, let $(S,\mc S)$ be a measure space and $(\xi)_{i \in I}$ a family of $S$-valued random variables on $\Omega$.  Let $A = L^{\infty}(\Omega)$, and let $\varphi:A \to \C$ be the expectation functional
\begin{equation*}
 \varphi(f) = \mb E[f].
\end{equation*}
Let $C$ be the algebra of bounded, complex-valued, $\mc S$-measurable functions on $S$.  For $i \in I$, define $\rho_i:C \to A$ by $\rho_i(f) = f \circ \xi_i$.  Then $\varphi_\rho$ is determined by
\begin{equation*}
 \varphi_\rho(f_1^{(i_1)}\dotsb f_k^{(i_k)}) = \mb E[ f_1(\xi_{i_1})\dotsb f_k(\xi_{i_k})]
\end{equation*}
for $f_1,\dotsc,f_k \in C$ and $i_1,\dotsc,i_k \in I$. 
\item   Let $C = \C[t]$, and let $(x_i)_{i \in I}$ be a family of self-adjoint random variables in $A$.  Define $\rho_i:C \to A$ to be the unique unital $*$-homomorphism such that $\rho_i(t) = x_i$.  Then $C_I = C \langle t_i: i \in I\rangle$, and we recover the usual definitions of the joint distribution and moments of the family $(x_i)_{i \in I}$.
 
\end{enumerate}
\end{remark*}

\begin{remark}
These definitions have natural ``operator-valued'' extensions given by replacing $\C$ by a more general algebra of scalars.  This is the right setting for the notion of freeness with amalgamation, which is the analogue of conditional independence in free probability.
\end{remark}

\begin{definition}
A \textit{$B$-valued probability space} $(A,E)$ consists of a unital $*$-algebra $A$, a $*$-subalgebra $1 \in B \subset A$, and a conditional expectation $E:A \to B$, i.e., $E$ is a linear map such that $E[1] = 1$ and
\begin{equation*}
 E[b_1ab_2] = b_1E[a]b_2
\end{equation*}
for all $b_1,b_2 \in B$ and $a \in A$. 
\end{definition}

\begin{definition}\label{bvalued}
Let $C$ be a unital $*$-algebra, $(A,E)$ a $B$-valued probability space and $(\rho_i)_{i \in I}$ a family of unital $*$-homomorphisms from $C$ into $A$.  
\begin{enumerate}
\item  We let $C_I^{B}$ denote the free product over $i \in I$, with amalgamation over $B$, of $C^{(i)} * B$, which is naturally isomorphic to $C_I * B$.  For each $i \in I$, we extend $\rho_i$ to a unital $*$-homomorphism $\widetilde \rho_i:C * B \to A$ by setting $\widetilde \rho_i = \rho_i * \mathrm{id}$. We then let $\widetilde \rho$ denote the induced unital $*$-homomorphism from $C_I^{B}$ into $A$, which is naturally identified with $\rho *\mathrm{id}$.  Explicitly, we have
\begin{equation*}
 \widetilde \rho(b_0c_1^{(i_1)}b_1\dotsb c_k^{(i_k)}b_k) = b_0\rho_{i_1}(c_1)b_1\dotsb \rho_{i_k}(c_k)b_k
\end{equation*}
for $b_0,\dotsc,b_k \in B$, $c_1,\dotsc,c_k \in C$ and $i_1,\dotsc,i_k \in I$.  
\item The \textit{$B$-valued joint distribution} of the family $(\rho_i)_{i \in I}$ is the linear map $E_\rho:C_I * B \to B$ defined by $E_\rho = E \circ \widetilde \rho$.  $E_\rho$ is determined by the \textit{$B$-valued moments}
\begin{equation*}
 E_\rho[b_0c_1^{(i_1)}\dotsb c_k^{(i_k)}b_k] = E[b_0\rho_{i_1}(c_1)\dotsb \rho_{i_k}(c_k)b_k]
\end{equation*}
for $c_1,\dotsc,c_k \in C$, $b_0,\dotsc,b_k \in B$ and $i_1,\dotsc,i_k \in I$.
\item The family $(\rho_i)_{i \in I}$ is called \textit{identically distributed with respect to $E$} if $E \circ \widetilde \rho_i = E \circ \widetilde \rho_j$ for all $i,j \in I$.  This is equivalent to the condition that
\begin{equation*}
 E[b_0\rho_i(c_1)\dotsb \rho_i(c_k)b_k] = E[b_0\rho_j(c_1)\dotsb \rho_j(c_k)b_k]
\end{equation*}
for any $i,j \in I$ and $c_1,\dotsc,c_k \in C$, $b_0,\dotsc,b_k \in B$.
\item The family $(\rho_i)_{i \in I}$ is called \textit{freely independent with respect to $E$}, or \textit{free with amalgamation over $B$}, if
\begin{equation*}
 E[\widetilde \rho_{i_1}(\beta_1)\dotsb \widetilde \rho_{i_k}(\beta_k)] = 0
\end{equation*}
whenever $i_1 \neq \dotsb \neq i_k \in I$, $\beta_1,\dotsc,\beta_k \in C * B$ and $E[\widetilde \rho_{i_l}(\beta_l)] = 0$ for $1 \leq l \leq k$.
\end{enumerate}
\end{definition}

\begin{remark}
Voiculescu introduced the notion of freeness with amalgamation and developed its basic theory in \cite{voi0}.  Freeness with amalgamation also has a rich combinatorial structure developed by Speicher \cite{sp1}.  The basic objects, which we will now recall, are non-crossing set partitions and free cumulants.  For further information on the combinatorial aspects of free probability, the reader is referred to the text \cite{ns}.  
\end{remark}

\begin{definition}\hfill

\begin{enumerate}
\renewcommand{\labelenumi}{(\roman{enumi})}
\item A \textit{partition} $\pi$ of a set $S$ is a collection of disjoint, non-empty sets $V_1,\dotsc,V_r$ such that $V_1 \cup \dotsb \cup V_r = S$.  $V_1,\dotsc,V_r$ are called the \textit{blocks} of $\pi$, and we set $|\pi| = r$. The collection of partitions of $S$ will be denoted $\mc P(S)$, or in the case that $S =\{1,\dotsc,k\}$ by $\mc P(k)$.
\item If $S$ is ordered, we say that $\pi \in \mc P(S)$ is \textit{non-crossing} if whenever $V,W$ are blocks of $\pi$ and $s_1 < t_1 < s_2 < t_2$ are such that $s_1,s_2 \in V$ and $t_1,t_2 \in W$, then $V = W$.  The set of non-crossing partitions of $S$ is denoted by $NC(S)$, or by $NC(k)$ in the case that $S = \{1,\dotsc,k\}$.

\item The non-crossing partitions can also be defined recursively, a partition $\pi \in \mc P(S)$ is non-crossing if and only if it has a block $V$ which is an interval, such that $\pi \setminus V$ is a non-crossing partition of $S \setminus V$.
\item Given $\pi,\sigma \in \mc P(S)$, we say that $\pi \leq \sigma$ if each block of $\pi$ is contained in a block of $\sigma$.  

\item  Given $i_1,\dotsc,i_k$ in some index set $I$, we denote by $\ker \mathbf i$ the element of $\mc P(k)$ whose blocks are the equivalence classes of the relation
\begin{equation*}
 s \sim t \Leftrightarrow i_s= i_t.
\end{equation*}
Note that if $\pi \in \mc P(k)$, then $\pi \leq \ker \mathbf i$ is equivalent to the condition that whenever $s$ and $t$ are in the same block of $\pi$, $i_s$ must equal $i_t$.

\end{enumerate}
\end{definition}

\begin{definition}
Let $(A,E)$ be a $B$-valued probability space. 
\begin{enumerate}
\renewcommand{\labelenumi}{(\roman{enumi})}
\item For each $k \in \N$, let $\rho^{(k)}:A^{\otimes_B k} \to B$ be a linear map (the tensor product is with respect to the natural $B-B$ bimodule structure on $A$).  For $n \in \N$ and $\pi \in NC(n)$, we define a linear map $\rho^{(\pi)}: A^{\otimes_B n} \to B$ recursively as follows.  If $\pi$ has only one block, we set
\begin{equation*}
 \rho^{(\pi)}[a_1 \otimes \dotsb \otimes a_n] = \rho^{(n)}(a_1 \otimes \dotsb \otimes a_n)
\end{equation*}
for any $a_1,\dotsc,a_n \in A$.  Otherwise, let $V = \{l+1,\dotsc,l+s\}$ be an interval of $\pi$.  We then define, for any $a_1,\dotsc,a_n \in A$, 
\begin{equation*}
 \rho^{(\pi)}[a_1 \otimes \dotsb \otimes a_n] = \rho^{(\pi \setminus V)}[a_1 \otimes \dotsb \otimes a_{l}\cdot \rho^{(s)}(a_{l+1} \otimes \dotsb \otimes a_{l+s}) \otimes \dotsb \otimes a_n].
\end{equation*}
For example, if
 \begin{equation*} 
\pi = \{\{1,5,8\},\{2,4\},\{3\},\{6,7\},\{9,10\}\} \in NC(10),
\end{equation*}
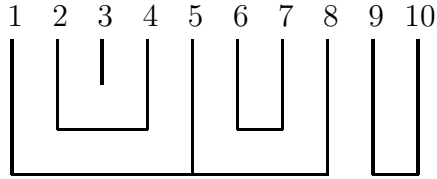
\begin{equation*}
 \setlength{\unitlength}{0.6cm} \begin{picture}(9,4)\thicklines \put(0,0){\line(0,1){3}}
\put(0,0){\line(1,0){7}} \put(8,0){\line(1,0){1}}\put(9,0){\line(0,1){3}} \put(8,0){\line(0,1){3}} \put(7,0){\line(0,1){3}} 
\put(1,1){\line(1,0){2}} \put(3,1){\line(0,1){2}}\put(1,1){\line(0,1){2}} \put(6,1){\line(0,1){2}}
\put(2,2){\line(0,1){1}} \put(3,2){\line(0,1){1}} \put(4,0){\line(0,1){3}}
\put(5,1){\line(0,1){2}} \put(5,1){\line(1,0){1}}
\put(-0.1,3.3){1} \put(0.9,3.3){2} \put(1.9,3.3){3}
\put(2.9,3.3){4} \put(3.9,3.3){5} \put(4.9,3.3){6} \put(5.9,3.3){7} \put(6.9,3.3){8}
\put(7.9,3.3){9} \put(8.7,3.3){10}
\end{picture}
\end{equation*}
then $\rho^{(\pi)}[a_1 \otimes \dotsb \otimes a_{10}]$ is given by
\begin{equation*}
 \rho^{(3)}(a_1\cdot \rho^{(2)}(a_2\cdot \rho^{(1)}(a_3) \otimes a_4) \otimes a_5 \cdot \rho^{(2)}(a_6 \otimes a_7) \otimes a_8)\cdot \rho^{(2)}(a_9 \otimes a_{10}).
\end{equation*}

\item For $k \in \N$, define the \textit{$B$-valued moment functions} $E^{(k)}:A^{\otimes_B k} \to B$ by
\begin{equation*}
 E^{(k)}[a_1 \otimes \dotsb \otimes a_k] = E[a_1\dotsb a_k].
\end{equation*}

\item The \textit{$B$-valued cumulant functions} $\kappa_E^{(k)}:A^{\otimes_B k} \to B$ are defined recursively for $\pi \in NC(k)$, $k \geq 1$, by the \textit{moment-cumulant formula}: for each $n \in \N$ and $a_1,\dotsc,a_n \in A$ we have
\begin{equation*}
 E[a_1\dotsb a_n] = \sum_{\pi \in NC(n)} \kappa_E^{(\pi)}[a_1 \otimes \dotsb \otimes a_n].
\end{equation*}
\end{enumerate} 
\end{definition}

\begin{remark*}
The cumulant functions can be solved for in terms of the moment functions by the following formula: for each $n \in \N$ and $a_1,\dotsc,a_n \in A$,
\begin{equation*}
 \kappa_E^{(\pi)}[a_1 \otimes \dotsb \otimes a_n] = \sum_{\substack{\sigma \in NC(n)\\ \sigma \leq \pi}} \mu_n(\sigma,\pi)E^{(\sigma)}[a_1 \otimes \dotsb \otimes a_n],
\end{equation*}
where $\mu_n$ is the \textit{M\"{o}bius function} on the partially ordered set $NC(n)$.  

The key relation between $B$-valued cumulant functions and free independence with amalgamation is that freeness can be characterized in terms of the ``vanishing of mixed cumulants''.
\end{remark*}

\begin{theorem}\textnormal{(\cite{sp1})} Let $C$ be a unital $*$-algebra, $(A,E)$ be a $B$-valued probability space and $(\rho_i)_{i \in I}$ a family of unital $*$-homomorphisms from $C$ into $A$.  Then the family $(\rho_i)_{i \in I}$ is free with amalgamation over $B$ if and only if
\begin{equation*}
 \kappa_{E}^{(\pi)}[\widetilde \rho_{i_1}(\beta_1)\otimes \dotsb \otimes \widetilde \rho_{i_k}(\beta_k)] = 0
\end{equation*}
whenever $i_1,\dotsc,i_k \in I$, $\beta_1,\dotsc,\beta_k \in C * B$ and $\pi \in NC(k)$ is such that $\pi \not\leq \ker \mathbf i$.
\end{theorem}

\begin{corollary}\label{vancum}
Let $C$ be a unital $*$-algebra, $(A,E)$ a $B$-valued probability space and $(\rho_i)_{i \in \N}$ a family of unital $*$-homomorphisms from $C$ into $A$.  Then $(\rho_i)_{i \in \N}$ is freely independent and identically distributed with respect to $E$ if and only if
\begin{equation*}
 E[\widetilde \rho_{i_1}(\beta_1)\dotsb \widetilde \rho_{i_k}(\beta_k)] = \sum_{\substack{\pi \in NC(k)\\ \pi \leq \ker \mathbf i}} \kappa_E^{(\pi)}[\widetilde \rho_{1}(\beta_1) \otimes \dotsb \otimes \widetilde \rho_{1}(\beta_k)]
\end{equation*}
for every $k \in \N$, $\beta_1,\dotsc,\beta_k \in C * B$ and $i_1,\dotsc,i_k \in I$. \qed
\end{corollary}

\begin{remark*}\textbf{Quantum Permutation Group.}  Wang introduced the following noncommutative analogue of $S_n$ in \cite{wang2}, and showed that it is the quantum automorphism group of a set with $n$ points.  For further information see \cite{bbc},\cite{bc2}.
\end{remark*}

\begin{definition}A matrix $(u_{ij})_{1 \leq i,j, \leq n} \in M_n(A)$, where $A$ is a unital C$^*$-algebra, is called a \textit{magic unitary} if

\begin{enumerate}
 \item $u_{ij}$ is a projection for each $1 \leq i,j \leq n$.
\item $u_{ik}u_{il} = 0$ and $u_{kj}u_{lj} = 0$ if $1 \leq i,j,k,l \leq n$ and $k \neq l$.
\item For each $1 \leq i,j \leq n$,
\begin{align*}
 \sum_{k=1}^n u_{ik} &= 1, & \sum_{k=1}^n u_{kj} &= 1.
\end{align*}
\end{enumerate}
Note that the second condition in fact follows from the third.  The \textit{quantum permutation group} $A_s(n)$ is defined as the universal C$^*$-algebra generated by elements $\{u_{ij}: 1 \leq i,j \leq n\}$ such that $(u_{ij})$ is a magic unitary.  $A_s(n)$ is a compact quantum group in the sense of Woronowicz \cite{wor1}, with comultiplication, counit and antipode given by
\begin{align*}
\Delta(u_{ij}) &= \sum_{k=1}^n u_{ik} \otimes u_{kj}\\
\epsilon(u_{ij}) &= \delta_{ij}\\
S(u_{ij}) &= u_{ji}.
\end{align*}
The existence of these maps is given by the universal property of $A_s(n)$.  
\end{definition}

\section{Quantum increasing sequences}\label{qinc}

In this section we introduce objects $A_i(k,n)$ which we call \textit{quantum increasing sequence spaces}.  As with Wang's quantum permutation group, the idea is to find a natural family of coordinates on the space of increasing sequences $1 \leq l_1 < \dotsb < l_k \leq n$ and ``remove commutativity''.

\begin{definition}
For $k, n \in \N$ with $k \leq n$, we define the \textit{quantum increasing sequence space} $A_i(k,n)$ to be the universal unital C$^*$-algebra generated by elements $\{u_{ij}: 1 \leq i \leq n, 1 \leq j \leq k\}$ such that
\begin{enumerate}
 \item $u_{ij}$ is an orthogonal projection: $u_{ij}^* = u_{ij} = u_{ij}^2$.
 \item each column of the rectangular matrix $u = (u_{ij})$ forms a partition of unity: for $1\leq j \leq k$ we have
\begin{equation*}
 \sum_{i = 1}^n u_{ij} = 1.
\end{equation*}
\item increasing sequence condition: 
\begin{equation*}
 u_{ij}u_{i'j'} = 0
\end{equation*}
if $j < j'$ and $i \geq i'$.
\end{enumerate}

\end{definition}

\begin{remark}
We note that the algebra $A_i(k,n)$, together with the morphism $\alpha:\C^n \to \C^k \otimes A_{i}(k,n)$ defined by
\begin{equation*}
 \alpha(e_i) = \sum_{j=1}^k e_j \otimes u_{ij},
\end{equation*}
gives a \textit{quantum family of maps} from $\{1,\dotsc,k\}$ to $\{1,\dotsc,n\}$, in the sense of So{\l}tan \cite{soltan}.
\end{remark}

The motivation for the above definition is as follows.  Consider the space $I_{k,n}$ of increasing sequences $\mathbf{l} = (1 \leq l_1 < \dotsb < l_k \leq n)$.  For $1 \leq i \leq n$, $1 \leq j \leq k$, define $f_{ij}:I_{k,n} \to \C$ by 
\begin{equation*}
 f_{ij}(\mathbf l) = \begin{cases} 1, &l_{j} = i\\0, &l_j \neq i\end{cases}.
\end{equation*}
The functions $f_{ij}$ generate $C(I_{k,n})$ by the Stone-Weierstrass theorem, and clearly satisfy the defining relations among the $u_{ij}$ above.  Moreover, it can be seen from the Gelfand theory that $C(I_{k,n})$ is the universal \textit{commutative} C$^*$-algebra generated by $\{f_{ij}: 1 \leq i \leq n, 1 \leq j \leq k\}$ satisfying these relations.  In other words, $C(I_{k,n})$ is the abelianization of $A_i(k,n)$.

\begin{remark}
A first question is whether $A_i(k,n)$ can be larger than $C(I_{k,n})$, i.e., ``do quantum increasing sequences exist''?  Clearly $A_i(k,n)$ is commutative and hence equal to $C(I_{k,n})$ for $k =1$.  Using Lemma \ref{zcoords} below, it is not hard to see that $A_i(k,n)$ is also commutative at $k = n$ and $n-1$.  In particular we have $A_i(k,n) = C(I_{k,n})$ whenever $n \leq 3$.  

However, if $p,q$ are arbitrary projections in any unital C$^*$-algebra then the following gives a representation of $A_i(2,4)$:
\begin{equation*}
 \begin{pmatrix}
  p & 0\\
1-p & 0\\
0 & q\\
0 & 1-q
 \end{pmatrix}
\end{equation*}
In particular, the free product $C(\Z_2) * C(\Z_2)$ is a quotient of $A_i(2,4)$ and hence $A_i(2,4)$ is infinite-dimensional.
\end{remark}

Observe that if $(1 \leq l_1 < \dotsb < l_k \leq n)$ then we must have $l_{j'} - l_{j} \geq j' - j$ for $1 \leq j \leq j' \leq k$.  In terms of the coordinates $f_{ij}$ on $C(I_{k,n})$, this means that $f_{ij}f_{i'j'} = 0$ if $i' -i < j' - j$.  This relation also holds for the coordinates $u_{ij}$ on $A_i(k,n)$, which will be useful to our further analysis.

\begin{lemma}\label{zcoords}
Fix $k,n \in \N$ with $k \leq n$, and let $\{u_{ij}:1 \leq i \leq n, 1 \leq j \leq k\}$ be the standard generators of $A_i(k,n)$.  Then 
\begin{enumerate}
\item $u_{ij}u_{i'j'} = 0$ if $1 \leq j \leq j' \leq k$ and $i' - i < j'-j$.  
 \item $u_{ij} = 0$ unless $j \leq i \leq n-k+j$, or equivalently $k+i - n \leq j \leq i$.
\end{enumerate}

\end{lemma}

\begin{proof}
(1) is trivial for $j=j'$, so fix $1 \leq j < j' \leq k$ and set $m = j'-j-1 \geq 0$.  Then we have
\begin{equation*}
 u_{ij}u_{i'j'} = u_{ij} \biggl(\prod_{l=1}^{m} \sum_{i_l=1}^n u_{i_l(j+l)}\biggr) u_{i'j'} = \sum_{1 \leq i_1,\dotsc,i_{m} \leq n} u_{ij}u_{i_1(j+1)}\dotsb u_{i_m(j+m)}u_{i'(j+m+1)}.
\end{equation*}
From the increasing sequence condition, each term in the sum is zero unless $i < i_1 < \dotsb < i_m < i'$, which implies $i' - i \geq m+1 = j'-j$.

For (2), note that from (1) we have $u_{l1}u_{ij} = 0$ if $i-l < j-1$, or equivalently $l > i-j + 1$.  So if $i < j$ then $u_{l1}u_{ij} = 0$ for $l=1,\dotsc,n$ and we then have
\begin{equation*}
 u_{ij} = \biggl(\sum_{l=1}^n u_{l1}\biggr) \cdot u_{ij} = 0.
\end{equation*}
Likewise we have $u_{ij}u_{lk} = 0$ if $l < k+i-j$, so if $i > n-k+j$ then this holds for $l =1,\dotsc,n$ and
\begin{equation*}
 u_{ij} = u_{ij} \cdot \biggl(\sum_{l=1}^n u_{lk}\biggr) = 0,
\end{equation*}
which completes the proof.
\end{proof}

Now observe that any increasing sequence $1 \leq l_1 < \dotsb < l_k \leq n$ can be extended to a permutation in $S_n$ which sends $j$ to $l_j$ for $1 \leq j \leq k$.  One way to create such an extension is to set $\pi(j) = l_j$ for $1 \leq j \leq k$, then inductively define $\pi(k+m)$, for $m=1,\dotsc,n-k$, by setting $\pi(k+m)$ to be the least element of $\{1,\dotsc,n\} \setminus\{\pi(1),\dotsc,\pi(k+m-1)\}$.  After a moment's thought, one sees that $m \leq \pi(k+m) \leq m+k$ and that $\pi(k+m) = m +p$ exactly when $l_{p} < m+p$ but $l_{p+1} > m+p$ for $1 \leq m \leq n-k$ and $0 \leq p \leq k$, where we set $l_0 = -\infty, l_{k+1} = \infty$.

This gives an inclusion of the space $I_{k,n}$ of increasing sequences into $S_n$, which dualizes to a unital $*$-homomorphism $C(S_n) \to C(I_{k,n})$.  Consider the natural coordinates $\{f_{ij}: 1\leq i,j \leq n\}$ on $S_n$ and $\{g_{ij}: 1 \leq i \leq n, 1 \leq j \leq k\}$ on $I_{k,n}$.  Clearly this map sends $f_{ij}$ to $g_{ij}$ for $1 \leq i \leq n$, $1 \leq j \leq k$.  From the remark at the end of the previous paragraph, it follows that $f_{i(k+m)}$ is sent to 0 unless $i = m+p$ for some $0 \leq p \leq k$, and that
\begin{equation*}
f_{(m+p)(k+m)} \mapsto \sum_{i=0}^{m+p-1} g_{ip} - g_{(i+1)(p+1)},
\end{equation*}
where we set $g_{00} = 1$ and $g_{i0} =g_{0i} = g_{i(k+1)} = 0$ for $i \geq 1$.  

For example, when $k=2$ and $n = 4$ the matrix $(f_{ij})$ is as follows:
\begin{equation*}
 \begin{pmatrix}
  g_{11} & 0 & 1-g_{11} & 0\\
g_{21} & g_{22} & g_{11}-g_{22} & 1-g_{11}-g_{21}\\
g_{31} & g_{32} & g_{22} & g_{11}+g_{21} - g_{22} - g_{32}\\
0 & g_{42} & 0 & g_{22}+g_{32}
 \end{pmatrix}
\end{equation*}

We can now use this formula to define a $*$-homomorphism from $A_s(n)$ to $A_i(k,n)$, which we might think of as ``extending quantum increasing sequences to quantum permutations''.

\begin{proposition}\label{qext}
Fix natural numbers $k < n$.  Let $\{v_{ij}: 1 \leq i \leq n, 1 \leq j \leq k\}$, $\{u_{ij}: 1 \leq i,j \leq n\}$ be the standard generators of $A_i(k,n)$, $A_s(n)$, respectively.  Then there is a unique unital $*$-homomorphism from $A_s(n)$ to $A_i(k,n)$ determined by
\begin{itemize}
 \item $u_{ij} \mapsto v_{ij}$ for $1 \leq i \leq n$, $1 \leq j \leq k$.
 \item $u_{i(k+m)} \mapsto 0$ for $1 \leq m \leq n-k$ and $i < m$ or $i > m+k$.
\item For $1 \leq m \leq n-k$ and $0 \leq p \leq k$, 
\begin{equation*}
 u_{(m+p)(k+m)} \mapsto \sum_{i=0}^{m+p-1} v_{ip} - v_{(i+1)(p+1)},
\end{equation*}
where we set $v_{00} = 1$ and $v_{i0} = v_{0i} = v_{i(k+1)} = 0$ for $i \geq 1$.
\end{itemize}

\end{proposition}

\begin{proof}
Let $(v_{ij})$ be the standard generators of $A_i(k,n)$, and define $\{u_{ij}: 1 \leq i,j \leq n\}$ in $A_i(k,n)$ by 
\begin{itemize}
\item $u_{ij} = v_{ij}$ for $1 \leq i \leq n$, $1 \leq j \leq k$.
 \item $u_{i(k+m)} = 0$ for $1 \leq m \leq n-k$ and $i < m$ or $i > m+k$.
\item For $1 \leq m \leq n-k$ and $0 \leq p \leq k$, 
\begin{equation*}
 u_{(m+p)(k+m)} = \sum_{i=0}^{m+p-1} v_{ip} - v_{(i+1)(p+1)},
\end{equation*}
where we set $v_{00} = 1$ and $v_{i0} = v_{0i} = v_{i(k+1)} = 0$ for $i \geq 1$.
\end{itemize}
We need to show that $(u_{ij})_{1 \leq i,j \leq n}$ satisfies the magic unitary condition, and the result will then follow from the universal property of $A_s(n)$.

First let us check that $u_{ij}$ is an orthogonal projection for $1 \leq i,j \leq n$.  The only non-trivial case is $u_{(m+p)(k+m)}$ for $1 \leq m \leq n-k$ and $0 \leq p \leq k$.  Here we just need to check that 
\begin{equation*}
 v_{l(p+1)} \leq \sum_{i=0}^{m+p-1} v_{ip}
\end{equation*}
for $1 \leq l \leq m+p$.  The cases $p = 0, k$ are trivial, so let $0 < p < k$.  We have
\begin{equation*}
 v_{l(p+1)} = v_{l(p+1)} \cdot \sum_{i=1}^n v_{ip} = v_{l(p+1)} \cdot \sum_{i=1}^{l-1} v_{ip},
\end{equation*}
where we have applied the increasing sequence condition $v_{l(p+1)}v_{ip} = 0$ for $i \geq l$.  So we have
\begin{equation*}
 v_{l(p+1)} \leq \sum_{i=1}^{l-1} v_{ip} \leq \sum_{i=0}^{m+p-1} v_{ip}
\end{equation*}
as desired.

Now we need to check that the sum along any row or column of $(u_{ij})$ gives the identity.  For the first $k$ columns, this follows from the defining relations of $v_{ij}$.  For $m= 1,\dotsc,n-k$, the sum along column $k+m$ gives
\begin{align*}
\sum_{l=1}^n u_{l(k+m)} &= \sum_{p=0}^{k} u_{(m+p)(k+m)} \\
&= \sum_{p=0}^k \sum_{i=0}^{m+p-1}v_{ip} - v_{(i+1)(p+1)}\\
\end{align*}
Now since $v_{ip} = v_{(i+1)(p+1)} = 0$ if $i < p$, we continue with
\begin{align*}
 \sum_{p=0}^k \sum_{i=p}^{m+p-1}v_{ip} - v_{(i+1)(p+1)} &= \sum_{p =0}^k \sum_{i=0}^{m-1} v_{(i+p)p} - v_{(i+p+1)(p+1)}\\
&= \sum_{i=0}^{m-1} \sum_{p=0}^k v_{(i+p)p} - v_{(i+p+1)(p+1)}\\
&= \sum_{i=0}^{m-1} v_{i0} - v_{(i+k+1)(k+1)}\\
&= 1,
\end{align*}
since the only nonzero term in the last sum is $v_{00} = 1$.

It now remains only to show that the sum along any row of $(u_{ij})$ gives the identity.  We have
\begin{align*}
 \sum_{j=1}^n u_{ij} &= \sum_{j=1}^k u_{ij} + \sum_{m = 1}^{n-k} \sum_{\substack{0 \leq p \leq k\\ m +p = i}} u_{i(k+m)}\\
&=\sum_{j=1}^k u_{ij} + \sum_{m= \max\{i-k,1\}}^{\min\{i,n-k\}} u_{i(k+m)}\\
&= \sum_{j=1}^k v_{ij} + \sum_{m=\max\{i-k,1\}}^{\min\{i,n-k\}} \sum_{l=0}^{i-1} v_{l(i-m)}-v_{(l+1)(i-m+1)}\\
&= \sum_{j=1}^k v_{ij} +  \biggl(\; \sum_{m=\max\{i-k,1\}}^{\min\{i,n-k\}} v_{0(i-m)}-v_{i(i-m+1)}\biggr) + \sum_{l=1}^{i-1} \sum_{m=\max\{i-k,1\}}^{\min\{i,n-k\}} v_{l(i-m)}-v_{l(i-m+1)}\\
&= \sum_{j=1}^k v_{ij} + \biggl(\; \sum_{m=\max\{i-k,1\}}^{\min\{i,n-k\}} v_{0(i-m)}-v_{i(i-m+1)}\biggr) + \sum_{l=1}^{i-1} v_{l\max\{0,k+i-n\}} - v_{l\min\{k+1,i\}}.
\end{align*}
Now note that if $1 \leq l \leq i-1$ then $v_{l\min\{k+1,i\}} = 0$, indeed this is true by definition if $\min\{k+1,i\} = k+1$, and if $\min\{k+1,i\} = i$ then $v_{li} = 0$ since $l < i$.  Also we have $v_{ij} = 0$ unless $k+i-n \leq j \leq i$.  Plugging this in above and rearranging terms, we have
\begin{equation*}
 \sum_{j = \max\{1,k+i-n\}}^{\min\{k,i\}}v_{ij} - \sum_{m=\max\{i-k,1\}}^{\min\{i,n-k\}} v_{i(i-m+1)} + \sum_{m=\max\{i-k,1\}}^{\min\{i,n-k\}} v_{0(i-m)} + \sum_{l=1}^{i-1} v_{l\max\{0,k+i-n\}}.
\end{equation*}
After reindexing the second sum and combining with the first, we obtain
\begin{equation*}
 \sum_{j=\max\{1,k+i-n\}}^{\max\{1,k+i+1-n\}-1} v_{ij} + \sum_{m=\max\{i-k,1\}}^{\min\{i,n-k\}} v_{0(i-m)} + \sum_{l=1}^{i-1} v_{l\max\{0,k+i-n\}}.
\end{equation*}
Now if $i \leq n-k$, then the first and third sums are zero while the second is 1.  If $i > n - k$ then the second sum is zero and the first and third combine as
\begin{equation*}
 \sum_{l=1}^{i}v_{l(k+i-n)}.
\end{equation*}
Now since $v_{l(k+i-n)} = 0$ if $l > n - k + (k+i-n) = i$, we have
\begin{equation*}
  \sum_{l=1}^{i}v_{l(k+i-n)} =  \sum_{l=1}^{n}v_{l(k+i-n)} = 1.
\end{equation*}
So $(u_{ij})$ does indeed satisfy the magic unitary condition, which completes the proof.
\end{proof}

\section{Quantum invariant sequences of random variables}\label{qinvariant}

In this section we introduce the notions of quantum exchangeability and quantum spreadability for sequences of noncommutative random variables, and prove the implications (i)$ \Rightarrow$ (ii)  and (iii) $\Rightarrow$ (i) in Theorem \ref{mainthm}.  First let us recall the notion of quantum exchangeability from \cite{ksp} (see also \cite{cur3}).

Let $C$ be a unital $*$-algebra.  For each $n \in \N$ there is a unique unital $*$-homomorphism $\alpha_n:C_n \to C_n \otimes A_s(n)$ determined by
\begin{equation*}
 \alpha_n(c^{(j)}) = \sum_{i=1}^n c^{(i)} \otimes u_{ij}
\end{equation*}
for $c \in C$ and $1 \leq j \leq n$, indeed this follows from the relations in $A_s(n)$ and the universal property of the free product $C_n = C^{(1)} * \dotsb * C^{(n)}$.  Moreover $\alpha_n$ is a right coaction of $A_s(n)$ in the sense that
\begin{align*}
 (\alpha_n \otimes \mathrm{id}) \circ \alpha_n &= (\mathrm{id} \otimes \alpha_n) \circ \alpha_n\\
(\mathrm{id} \otimes \epsilon) \circ \alpha_n &= \mathrm{id},
\end{align*}
see \cite{cur3} for details.  The coaction $\alpha_n$ may be regarded as ``quantum permuting'' the $n$ copies of $C$ inside $C_n$.

\begin{definition}
Let $C$ be a unital $*$-algebra, $(A,\varphi)$ a noncommutative probability space and $(\rho_1,\dotsc,\rho_n)$ a sequence of unital $*$-homomorphisms of $C$ into $A$.  We say that the distribution $\varphi_\rho$ is \textit{invariant under quantum permutations}, or that the sequence is \textit{quantum exchangeable}, if $\varphi_\rho$ is invariant under the coaction $\alpha_n$, i.e.,
\begin{equation*}
 (\varphi_\rho \otimes \mathrm{id})\alpha_n(c) = \varphi_\rho(c)1_{A_s(n)}
\end{equation*}
for any $c \in C_n$.  

This is extended to infinite sequences $(\rho_i)_{i \in \N}$ by requiring that $(\rho_1,\dotsc,\rho_n)$ is quantum exchangeable for each $n \in \N$.
\end{definition}

\begin{remark*}\textit{Remarks}.\label{qexcdef}
\begin{enumerate}
\item More explicitly, this amounts to the condition that
\begin{equation*}
 \sum_{1 \leq i_1,\dotsc,i_k \leq n} \varphi(\rho_{i_1}(c_1)\dotsb \rho_{i_k}(c_k))u_{i_1j_1}\dotsb u_{i_kj_k} = \varphi(\rho_{j_1}(c_1)\dotsb \rho_{j_k}(c_k))\cdot 1
\end{equation*}
for any $c_1,\dotsc,c_k \in C$ and $1 \leq j_1,\dotsc,j_k \leq n$, where $u_{ij}$ are the standard generators of $A_s(n)$.
\item By the universal property of $A_s(n)$, the sequence $(\rho_1,\dotsc,\rho_n)$ is quantum exchangeable if and only if the equation in (1) holds for any family $\{u_{ij}:1 \leq i,j \leq n\}$ of projections in a unital C$^*$-algebra $B$ such that $(u_{ij}) \in M_n(B)$ is a magic unitary matrix.
\item For $1 \leq i,j \leq n$, define $f_{ij} \in C(S_n)$ by $f_{ij}(\pi) = \delta_{i\pi(j)}$.  The matrix $(f_{ij})$ is a magic unitary, and the equation in (1) becomes
\begin{align*}
 \varphi(\rho_{j_1}(c_1)\dotsb \rho_{j_k}(c_k))1_{C(S_n)} &= \sum_{1 \leq i_1,\dotsc,i_k \leq n} \varphi(\rho_{i_1}(c_1)\dotsb \rho_{i_n}(c_n))f_{i_1j_1}\dotsb f_{i_kj_k}.
\end{align*}
Evaluating both sides at $\pi \in S_n$, we find
\begin{equation*}
\varphi(\rho_{j_1}(c_1)\dotsb \rho_{j_k}(c_k)) = \varphi(\rho_{\pi(j_1)}(c_1)\dotsb \rho_{\pi(j_k)}(c_k)),
\end{equation*}
so that quantum exchangeability implies invariance under classical permutations.
\end{enumerate}

\end{remark*}

It is shown in \cite{ksp} that any sequence $(\rho_1,\dotsc,\rho_n)$ which is freely independent and identically distributed with respect to a conditional expectation which preserves $\varphi$ is quantum exchangeable.  For the convenience of the reader we include a sketch of the proof, and refer to that paper for details.  Note that the implication (iii) $\Rightarrow$ (i) in Theorem \ref{mainthm} follows immediately.

\begin{proposition}
Let $C$ be a unital $*$-algebra and $(\rho_1,\rho_2,\dotsc,\rho_n)$ a sequence of unital $*$-homomorphisms from $C$ into a noncommutative probability space $(A,\varphi)$.  Let $B \subset A$ be a unital $*$-subalgebra and suppose that there is a $\varphi$-preserving conditional expectation $E:A \to B$ such that $(\rho_1,\dotsc,\rho_n)$ is freely independent and identically distributed with respect to $E$.  Then $(\rho_1,\dotsc,\rho_n)$ is quantum exchangeable.
\end{proposition}

\begin{proof}
Let $c_1,\dotsc,c_k \in C$ and $1 \leq j_1,\dotsc,j_k \leq n$.  We have
\begin{multline*}
 \sum_{1 \leq i_1,\dotsc,i_k \leq n} \negthickspace \varphi(\rho_{i_1}(c_1)\dotsb \rho_{i_k}(c_k))u_{i_1j_1}\dotsb u_{i_kj_k} = \negthickspace \sum_{1 \leq i_1,\dotsc,i_k \leq n} \varphi(E[\rho_{i_1}(c_1)\dotsb \rho_{i_k}(c_k)])u_{i_1j_1}\dotsb u_{i_kj_k}\\
= \sum_{1 \leq i_1,\dotsc,i_k \leq n} \sum_{\substack{\pi \in NC(k)\\ \pi \leq \ker \mathbf i}} \varphi(\kappa_E^{(\pi)}[\rho_{1}(c_1) \otimes \dotsb \otimes \rho_{1}(c_k)])u_{i_1j_1}\dotsb u_{i_kj_k}\\
= \sum_{\pi \in NC(k)} \varphi(\kappa_E^{(\pi)}[\rho_1(c_1) \otimes \dotsb \otimes \rho_1(c_k)] )\sum_{\substack{1 \leq i_1,\dotsc,i_k \leq n\\ \pi \leq \ker \mathbf i}} u_{i_1j_1}\dotsb u_{i_kj_k},
\end{multline*}
where in the second line we have applied Corollary \ref{vancum}.  It can be seen from induction on the number of blocks of $\pi$ that
\begin{equation*}
 \sum_{\substack{1 \leq i_1,\dotsc,i_k \leq n\\ \pi \leq \ker \mathbf i}} u_{i_1j_1}\dotsb u_{i_kj_k} = \begin{cases} 1_{A_s(n)}, & \pi \leq \ker \mathbf j\\ 0, & \text{otherwise}\end{cases},
\end{equation*}
and it follows that
\begin{align*}
 \sum_{1 \leq i_1,\dotsc,i_k \leq n} \varphi(\rho_{i_1}(c_1)\dotsb \rho_{i_k}(c_k))u_{i_1j_1}\dotsb u_{i_kj_k} &= \sum_{\substack{\pi \in NC(k)\\ \pi \leq \ker \mathbf j}} \varphi(\kappa_E^{(\pi)}[\rho_1(c_1) \otimes \dotsb \otimes \rho_1(c_k)] )1_{A_s(n)}\\
&= \varphi(\rho_{j_1}(c_1)\dotsb \rho_{j_k}(c_k))1_{A_s(n)},
\end{align*}
where again we have applied Corollary \ref{vancum}.
\end{proof}

We will now introduce the quantum spreadability condition.  Let $C$ be a unital $*$-algebra, then for any natural numbers $k \leq n$ there is a unique unital $*$-homomorphism $\alpha_{k,n}:C_k \to C_n \otimes A_{i}(k,n)$ determined by
\begin{equation*}
 \alpha_{k,n}(c^{(j)}) = \sum_{i=1}^n c^{(i)} \otimes u_{ij}
\end{equation*}
for $c \in C$ and $1 \leq j \leq k$, indeed this follows as above from the relations in $A_i(k,n)$ and the universal property of $C_k$.

\begin{definition}
Let $C$ be a unital $*$-algebra and $(\rho_1,\rho_2,\dotsc,\rho_n)$ a sequence of unital $*$-homomorphisms from $C$ into a noncommutative probability space $(A,\varphi)$.  We say that the distribution $\varphi_\rho$ is \textit{invariant under quantum spreading}, or that the sequence is \textit{quantum spreadable}, if for each $k = 1,\dotsc,n$ the distribution $\varphi_\rho$ is invariant under $\alpha_{k,n}$ in the sense that
\begin{equation*}
 (\varphi_\rho \otimes \mathrm{id})\alpha_{k,n}(c) = \varphi_\rho(c)1_{A_{i}(k,n)}
\end{equation*}
for any $c \in C_k$. 

An infinite sequence $(\rho_1,\rho_2,\dotsc)$ is called \textit{quantum spreadable} if $(\rho_1,\dotsc,\rho_n)$ is quantum spreadable for each $n$.
\end{definition}

\begin{remark}\hfill
\begin{enumerate}
\item Explicitly, the condition is that for each $k = 1,\dotsc,n$ we have
\begin{equation*}\label{qspread}
 \varphi(\rho_{j_1}(c_1)\dotsb \rho_{j_m}(c_m)) \cdot 1 = \sum_{1 \leq i_1,\dotsc,i_m \leq n} \varphi(\rho_{i_1}(c_1)\dotsb \rho_{i_m}(c_m)) \cdot u_{i_1j_1}\dotsb u_{i_mj_m}
\end{equation*}
for all $1 \leq j_1,\dotsc,j_m \leq k$ and $c_1,\dotsc,c_m \in C$, where $(u_{ij})$ denote the standard generators of $A_i(k,n)$.  

 \item From the universal property of $A_i(k,n)$, the sequence $(\rho_1,\dotsc,\rho_n)$ is quantum spreadable if and only if for each $1\leq k \leq n$, equation (\ref{qspread}) holds for any family $\{u_{ij}:1 \leq i \leq n, 1 \leq j \leq k\}$ of projections in a unital C$^*$-algebra $B$ which satisfy the definining relations of $A_i(k,n)$.
 \item Let $(f_{ij})$ denote the generators of $C(I_{k,n})$ introduced in Section \ref{qinc}.  Plugging $f_{ij}$ into equation (\ref{qspread}) and applying both sides to $\mathbf l = (1 \leq l_1 < \dotsb < l_k \leq n)$, we have
\begin{align*}
 \varphi(\rho_{j_1}(c_1)\dotsb \rho_{j_m}(c_m)) &= \sum_{1 \leq i_1,\dotsc,i_m \leq n} \varphi(\rho_{i_1}(c_1)\dotsb \rho_{i_m}(c_m)f_{i_1j_1}(\mathbf l)\dotsb f_{i_mj_m}(\mathbf l)\\
&= \varphi(\rho_{l_{j_1}}(c_1)\dotsb \rho_{l_{j_m}}(c_m))
\end{align*}
for any $1 \leq j_1,\dotsc,j_m \leq k$.  So $(\rho_1,\dotsc,\rho_k)$ has the same distribution as $(\rho_{l_1},\dotsc,\rho_{l_k})$, and hence quantum spreadability implies classical spreadability.  In particular, quantum spreadable sequences are identically distributed.
\end{enumerate}

\end{remark}
We can now prove the implication (i) $\Rightarrow$ (ii) of Theorem \ref{mainthm}, this holds in fact for finite sequences and in a purely algebraic context:

\begin{proposition}
Let $C$ be a unital $*$-algebra and $(\rho_1,\dotsc,\rho_n)$ be a sequence of unital $*$-homomorphisms from $C$ into a noncommutative probability space $(A,\varphi)$.  If the sequence $(\rho_1,\dotsc,\rho_n)$ is quantum exchangeable, then it is quantum spreadable.
\end{proposition}

\begin{proof}
Fix $1 \leq k \leq n$ and let $\{v_{ij}: 1\leq i \leq n, 1 \leq j \leq k\}$ and $\{u_{ij}: 1\leq i,j \leq n\}$ be the standard generators of $A_i(k,n)$ and $A_s(n)$, respectively.  Assume $(\rho_1,\dotsc,\rho_n)$ is quantum exchangeable, and fix $1 \leq j_1,\dotsc,j_m \leq k$ and $c_1,\dotsc,c_m \in C$.  We have
\begin{equation*}
 \varphi(\rho_{j_1}(c_1)\dotsb \rho_{j_m}(c_m))1_{A_s(n)} = \sum_{1 \leq i_1,\dotsc,i_m \leq n} \varphi(\rho_{i_1}(c_1)\dotsb \rho_{i_m}(c_m)) \cdot u_{i_1j_1}\dotsb u_{i_mj_m}.
\end{equation*}
By Proposition \ref{qext}, there is a unital $*$-homomorphism from $A_s(n)$ to $A_i(k,n)$ which sends $u_{ij}$ to $v_{ij}$ for $1 \leq i \leq n$, $1 \leq j \leq k$.  Applying this map to both sides of the above equation, we obtain
\begin{equation*}
  \varphi(\rho_{j_1}(c_1)\dotsb \rho_{j_m}(c_m))1_{A_i(k,n)} = \sum_{1 \leq i_1,\dotsc,i_m \leq n} \varphi(\rho_{i_1}(c_1)\dotsb \rho_{i_m}(c_m)) \cdot v_{i_1j_1}\dotsb v_{i_mj_m},
\end{equation*}
so that $(\rho_1,\dotsc,\rho_n)$ is quantum spreadable as desired.
\end{proof}

\section{Quantum spreadability implies freeness with amalgamation}\label{ryll}

\begin{remark*}
In this section we will complete the proof of Theorem \ref{mainthm}.  Throughout this section we will assume that $C$ is a unital $*$-algebra, and that $(\rho_i)_{i \in \N}$ is an infinite sequence of unital $*$-homomorphisms from $C$ into a W$^*$-probability space $(M,\tau)$.  $B$ will denote the tail algebra:
\begin{equation*}
 B = \bigcap_{n \geq 1} W^*\bigl(\{\rho_i(c):c \in C, i \geq n\}\bigr).
\end{equation*}
$L^2(M)$ will denote the Hilbert space given by the GNS-representation for $\tau$.  Since $\tau$ is a trace, there is a unique conditional expectation $E:M \to B$ given my $E[m] = P(m)$, where $P$ is the orthogonal projection of $L^2(M)$ onto $L^2(B)$.  

We will assume without loss of generality that $M$ is generated by $\rho_\infty(C_\infty)$, i.e.,
\begin{equation*}
 M = W^*\bigl( \{\rho_i(c): i \in I, c \in C\} \bigr).
\end{equation*}
Observe that if the sequence $(\rho_i)_{i \in \N}$ is spreadable and hence stationary,  the linear map determined by
\begin{equation*}
 U(\rho_{i_1}(c_1)\dotsb \rho_{i_m}(c_m)) = \rho_{i_1+1}(c_1)\dotsb \rho_{i_m+1}(c_m)
\end{equation*}
for $i_1,\dotsc,i_m \in \N$ and $c_1,\dotsc,c_m \in C$, is well-defined and extends to an isometry $U:L^2(M) \to L^2(M)$.  
\end{remark*}

Recall from Definition \ref{bvalued} that we set $\widetilde \rho_i = \rho_i * \mathrm{id}:C *B \to M$.  We will begin by showing that if $(\rho_i)_{i \in \N}$ is quantum spreadable, then the $B$-valued distribution of $(\widetilde \rho_i)_{i \in \N}$ is also invariant under quantum spreading.  By this we mean that the joint distribution $E_\rho$ is invariant under the $*$-homomorphisms $\widetilde \alpha_{k,n}: C_k * B \to (C_n * B) \otimes A_{i}(k,n)$ determined by
\begin{equation*}
 \widetilde \alpha_{k,n} (b_0c_1^{(j_1)}b_1 \dotsb c_m^{(j_m)}b_m) = \sum_{1 \leq i_1,\dotsc,i_m \leq n} b_0c_1^{(i_1)}b_1\dotsb c_m^{(i_m)}b_m \otimes u_{i_1j_1}\dotsb u_{i_mj_m}
\end{equation*}
for all $k \leq n$, $1 \leq j_1,\dotsc,j_m \leq k$, $b_0,\dotsc,b_m \in B$ and $c_1,\dotsc,c_m \in C$.

Note that if $1 \leq j \leq k$, $b_0,\dotsc,b_m \in B$ and $c_1,\dotsc,c_m \in C$ then
\begin{align*}
 \widetilde \alpha_{k,n}(b_0c_1^{(j)}\dotsb c_m^{(j)}b_m) &= \sum_{1 \leq i_1,\dotsc,i_m \leq n} b_0c_1^{(i_1)}\dotsb c_m^{(i_m)}b_m \otimes u_{i_1j}\dotsb u_{i_mj}\\
&= \sum_{i=1}^n b_0c_1^{(i)}\dotsb c_m^{(i)}b_m \otimes u_{ij},
\end{align*}
from which it follows that if $\beta \in C * B$ then
\begin{equation*}
 \widetilde \alpha_{k,n}(\beta^{(j)}) = \sum_{i=1}^{n} \beta^{(i)} \otimes u_{ij}.
\end{equation*}

\begin{proposition}\label{expspread}
Suppose that the sequence $(\rho_i)_{i \in \N}$ is quantum spreadable.  Then the joint distribution of $(\widetilde{\rho}_i)_{i \in \N}$ with respect to $E$ is invariant under quantum spreading.  Explicitly, for each $k \leq n$, $1 \leq j_1,\dotsc,j_m \leq k$ and $\beta_1,\dotsc,\beta_m \in C * B$ we have
\begin{equation*}
 E[\widetilde{\rho}_{j_1}(\beta_1)\dotsb \widetilde{\rho}_{j_m}(\beta_m)] \otimes 1_{A_i(k,n)} = \sum_{1 \leq i_1,\dotsc,i_m \leq n} E[\widetilde{\rho}_{i_1}(\beta_1)\dotsb \widetilde{\rho}_{i_m}(\beta_m)] \otimes u_{i_1j_1}\dotsb u_{i_mj_m},
\end{equation*}
where the equality holds in $B \otimes A_i(k,n)$.
\end{proposition}

\begin{proof}
We need to show that if $1 \leq j_1,\dotsc,j_m \leq k$, $b_0,\dotsc,b_m \in B$ and $c_1,\dotsc,c_m \in C$ then
\begin{equation*}
 E[b_0\rho_{j_1}(c_1)\dotsb \rho_{j_m}(c_m)b_m] \otimes 1 = \sum_{1 \leq i_1,\dotsc,i_m \leq n} E[b_0\rho_{i_1}(c_1)\dotsb \rho_{i_m}(c_m)b_m] \otimes u_{i_1j_1}\dotsb u_{i_mj_m}.
\end{equation*}
Since $E$ preserves the faithful state $\tau$, it suffices to show that
\begin{equation*}
 \tau(b_0\rho_{j_1}(c_1)\dotsb \rho_{j_m}(c_m)b_m) \otimes 1 = \sum_{1 \leq i_1,\dotsc,i_m \leq n} \tau(b_0\rho_{i_1}(c_1)\dotsb \rho_{i_m}(c_m)b_m) \otimes u_{i_1j_1}\dotsb u_{i_mj_m}.
\end{equation*}
We will show that this in fact holds for $b_0,\dotsc,b_m$ in $W^*(\{\rho_i(c): i > k, c \in C\})$.  By Kaplansky's density theorem, it suffices to consider the case that $b_0,\dotsc,b_m$ are elements of the form $\rho_{l_1}(d_1)\dotsb \rho_{l_r}(d_r)$ for $k < l_1,\dotsc,l_r \leq N$ and $d_1,\dotsc,d_r \in C$.  

To show this, we extend $(u_{ij})$ to a $(n+N) \times (k+N)$ matrix by setting
\begin{equation*}
 v_{ij} = \begin{cases} u_{ij}, & 1 \leq i \leq n, 1 \leq j \leq k\\
           \delta_{(i-n)(j-k)}, &i > n, j > k\\
	    0, & \text{otherwise}
          \end{cases}
\end{equation*}
Observe if $b = \rho_{l_1}(d_1)\dotsc \rho_{l_r}(d_r)$ is as above, then
\begin{align*}
 \sum_{1 \leq i_1,\dotsc,i_r \leq n+N} \rho_{i_1}(d_1)\dotsb \rho_{i_m}(d_r) \otimes v_{i_1l_1}\dotsb v_{i_rl_r} &= \rho_{l_1+(n-k)}(d_1)\dotsb \rho_{l_r+(n-k)}(d_r) \otimes 1_{A_i(k,n)}\\
&= U^{(n-k)}(b) \otimes 1_{A_i(k,n)}.
\end{align*}
Now it is clear that $(v_{ij})$ satisfies the defining relations of $A_{i}(k+N,n+N)$, so applying the quantum spreadability condition with $(v_{ij})$, we have
\begin{multline*}
  \tau(b_0\rho_{j_1}(c_1)\dotsb \rho_{j_m}(c_m)b_m)\\
 =  \sum_{1 \leq i_1,\dotsc,i_m \leq n} \tau\bigl(U^{(n-k)}(b_0)\rho_{i_1}(c_1)\dotsb \rho_{i_m}(c_m)U^{(n-k)}(b_m)\bigr) \otimes u_{i_1j_1}\dotsb u_{i_mj_m}.
\end{multline*}
But since $(\rho_i)_{i \in \N}$ is spreadable, the right hand side is equal to
\begin{equation*}
 \sum_{1 \leq i_1,\dotsc,i_m \leq n} \tau\bigl(b_0\rho_{i_1}(c_1)\dotsb \rho_{i_m}(c_m)b_m\bigr) \otimes u_{i_1j_1}\dotsb u_{i_mj_m},
\end{equation*}
which completes the proof.

\end{proof}

\begin{remark*}
The key ingredient in our proof that an infinite quantum spreadable sequence is free with amalgamation is a ``measure'' on the space of quantum increasing sequences, i.e., a state on $A_i(k,n)$.  Unlike in the classical case, there does not appear to be a good notion of ``uniform'' measure on this quantum space.  Instead, we will use the measures induced by a certain representation of $A_i(k,k\cdot n)$.  
\end{remark*}

\begin{proposition}\label{measure}
Fix $k,n \in \N$.  Then there is a state $\psi_{k,n}:A_i(k,k\cdot n) \to \C$ such that:
\begin{enumerate}
\item
\begin{equation*}
 \psi_{k,n}(u_{l_1j_1}\dotsb u_{l_mj_m}) = 0
\end{equation*}
unless $(j_r-1)\cdot n < l_r \leq j_r\cdot n$ for $r=1,\dotsc,m$.
\item
\begin{equation*}
 \psi_{k,n}(u_{((j_1-1)\cdot n + i_1)j_1}\dotsb u_{((j_m-1)\cdot n+i_m)j_m}) = \sum_{\substack{\pi \in NC(m)\\ \pi \leq \ker \mathbf j}} \sum_{\substack{\sigma \in NC(m) \\ \sigma \leq \pi \wedge \ker \mathbf i}} \mu_m(\sigma,\pi) n^{-|\sigma|}
\end{equation*}
for all $1 \leq j_1,\dotsc,j_m \leq k$ and $1\leq i_1,\dotsc,i_m \leq n$.
\end{enumerate}
\end{proposition}

\begin{proof}
Let $\{p_{ij}: 1 \leq i \leq n, 1 \leq j \leq k\}$ be projections in a C$^*$-probability space $(A,\varphi)$ such that 
\begin{enumerate}
\item The families $(\{p_{i1}: 1 \leq i \leq n\},\dotsc, \{p_{ik}:1 \leq i \leq n\})$ are freely independent.
\item For $j =1,\dotsc,k$, we have
\begin{equation*}
 \sum_{i=1}^n p_{ij} = 1,
\end{equation*}
and $\varphi(p_{ij}) = n^{-1}$ for $1 \leq i \leq n$.
\end{enumerate}

Define $\{u_{lj}: 1 \leq l \leq kn, 1 \leq j \leq k\}$ by $u_{lj} = 0$ unless $(j-1)\cdot n < l \leq j\cdot n$, and
\begin{equation*}
 u_{((j-1)\cdot n + i)j} = p_{ij}
\end{equation*}
for $1 \leq i \leq n$, so that $(u_{lj})$ is given by the following matrix:
\begin{equation*}
 \begin{pmatrix}
  p_{11} & 0 & \cdots & 0\\
\vdots & \vdots & \ddots & \vdots\\
p_{1n} & 0 & \cdots & 0\\
0 & p_{21} & \cdots & 0\\
\vdots & \vdots & \ddots & \vdots\\
0 & p_{2n} & \cdots & 0\\
0 & 0 & \cdots & 0\\
\vdots & \vdots & \ddots & \vdots\\
0 & 0 & \cdots & p_{k1}\\
\vdots & \vdots & \ddots & \vdots\\
0 & 0 & \cdots & p_{kn}
 \end{pmatrix}
\end{equation*}

  Clearly $(u_{lj})$ satisfies the defining relations of $A_i(k,k\cdot n)$ and so we obtain a unital $*$-homorphism from $A_i(k,k\cdot n)$ into $A$.  Composing with $\varphi$ gives a state $\psi_{k,n}:A_i(k,k\cdot n) \to \C$, and we need only show that $(u_{lj})$ in $(A,\varphi)$ has the distribution appearing in the statement.

(1) is trivial, as $u_{l_1j_1}\dotsb u_{l_mj_m} = 0$ unless $(j_r-1)\cdot n < l_r \leq j_r\cdot n$ for $r=1,\dotsc,m$.  For (2), we need to show that
\begin{equation*}
 \varphi(p_{i_1j_1}\dotsb p_{i_mj_m}) = \sum_{\substack{\pi \in NC(m)\\ \pi \leq \ker \mathbf j}} \sum_{\substack{\sigma \in NC(m) \\ \sigma \leq \pi \wedge \ker \mathbf i}} \mu_m(\sigma,\pi) n^{-|\sigma|}.
\end{equation*}
Now by freeness, we have
\begin{align*}
 \varphi(p_{i_1j_1}\dotsb p_{i_mj_m}) &= \sum_{\substack{\pi \in NC(m)\\ \pi \leq \ker \mathbf j}} \kappa^{(\pi)}[p_{i_1j_1} \otimes \dotsb \otimes p_{i_mj_m}]\\
&= \sum_{\substack{\pi \in NC(m)\\ \pi \leq \ker \mathbf j}} \sum_{\substack{\sigma \in NC(m)\\ \sigma \leq \pi}} \mu_m(\sigma,\pi)\varphi^{(\sigma)}[p_{i_1j_1} \otimes \dotsb \otimes p_{i_mj_m}].
\end{align*}
Now since for $1 \leq j, l_1,\dotsc,l_s \leq n$ we have 
\begin{equation*}
\varphi(p_{l_1j}\dotsb p_{l_sj}) = \begin{cases} n^{-1}, &l_1 = \dotsb = l_s\\ 0, & \text{otherwise} \end{cases},
\end{equation*}
it follows that if $\sigma \leq \ker \mathbf j$ then
\begin{equation*}
 \varphi^{(\sigma)}[p_{i_1j_1} \otimes \dotsb \otimes p_{i_mj_m}] = \begin{cases} n^{-|\sigma|}, &\sigma \leq \ker \mathbf i\\ 0, & \sigma \not\leq \ker \mathbf i \end{cases}.
\end{equation*}
Combining this with the previous equation yields the desired result.

\end{proof}

\begin{remark}
Observe that the formula in (2) above has a very similar structure to the highest order expansion of the \textit{Weingarten formula} for evaluating integrals over the quantum permutation group $A_s(n)$ with respect to its Haar state, see \cite{bc2,cur3}.  

The final tool which we require to complete the proof of Theorem \ref{mainthm} is von Neumann's mean ergodic theorem.  This will allow us to give a formula for the expectation functionals $E^{(\sigma)}$ as certain weighted averages.  We note that the unpleasant indices which appear are chosen as to correspond to the formula in Proposition \ref{measure}.  
\end{remark}

\begin{lemma}\label{ergodic}
Suppose that the sequence $(\rho_i)_{i \in \N}$ is quantum spreadable.  Then for any $j \in \N$ and $\beta \in C * B$, we have
\begin{equation*}
 E[\widetilde \rho_1(\beta)] = \lim_{n \to \infty} \frac{1}{n} \sum_{i=1}^n \widetilde \rho_{(j-1)\cdot n + i}(\beta),
\end{equation*}
with convergence in $|\;|_2$.
\end{lemma}

\begin{proof}
Since $(\rho_i)_{i \in \N}$ is spreadable, we have
\begin{equation*}
 \tau(m_1m_2) = \tau(m_1U(m_2))
\end{equation*}
whenever $m_1 \in W^*(\{\rho_i(c): 1 \leq i \leq n, c \in C\})$ and $m_2 \in W^*(\{\rho_{i}(c):i > n, c \in C\})$.  It follows that
\begin{equation*}
 \tau(mb) = \tau(mU(b))
\end{equation*}
for $m \in M$ and $b \in B$, hence $b = U(b)$.  It follows easily that
\begin{equation*}
 U(\widetilde \rho_i(\beta)) = \widetilde \rho_{i+1}(\beta)
\end{equation*}
for any $i \in \N$ and $\beta \in C * B$.  

Since it is clear that any vector fixed by $U$ must lie in $L^2(B)$, we have in fact the equality
\begin{equation*}
 L^2(B) = \{\xi \in L^2(M): U\xi = \xi\}.
\end{equation*}
By von Neumann's mean ergodic theorem, we have
\begin{equation*}
 \lim_{n \to \infty} \frac{1}{n}\sum_{i=0}^{n-1} U^i = P,
\end{equation*}
where $P$ is the orthogonal projection of $L^2(M)$ onto $L^2(B)$ and the limit holds in the strong operator topology.  Therefore for any $m \in M$ we have
\begin{equation*}
 E[m] = P(m) = \lim_{n \to \infty} \frac{1}{n}\sum_{i=0}^{n-1} U^i(m),
\end{equation*}
with the limit holding in $|\;|_2$.  Since $U$ is contractive in $|\;|_2$, we have also for any $j \in \N$ that
\begin{align*}
 \lim_{n \to \infty} \frac{1}{n} \sum_{i=0}^{n-1} U^{(j-1)\cdot n + i}(m) &= \lim_{n \to \infty} U^{(j-1)\cdot n} \biggl(\frac{1}{n}\sum_{i=0}^{n-1} U^i(m)\biggr)\\
&= \lim_{n \to \infty} U^{(j-1)\cdot n} P(m)\\
&= E[m],
\end{align*}
since $U \cdot P = P$.  Applying this to $m = \widetilde \rho_1(\beta)$ gives the desired result.
\end{proof}

\begin{proposition}\label{expform}
Suppose that the sequence $(\rho_i)_{i \in \N}$ is quantum spreadable.  Fix $j_1,\dotsc,j_m \in \N$ and choose $\sigma \in NC(m)$  such that $\sigma \leq \ker \mathbf j$.  Then for any $\beta_1,\dotsc,\beta_m \in C * B$, we have
\begin{equation*}
 E^{(\sigma)}[\widetilde{\rho}_1(\beta_1) \otimes \dotsb \otimes \widetilde{\rho}_1(\beta_m)] = \lim_{n \to \infty} n^{-|\sigma|} \sum_{\substack{1 \leq i_1,\dotsc,i_m \leq n\\ \sigma \leq \ker \mathbf i}} \widetilde{\rho}_{(j_1-1)\cdot n + i_1}(\beta_1)\dotsb \widetilde{\rho}_{(j_m -1)\cdot n + i_m}(\beta_m),
\end{equation*}
with convergence in $|\;|_2$.
\end{proposition}

\begin{proof}
We will use induction on the number of blocks of $\sigma$.  If $\sigma = 1_m$ has only one block, then $\sigma \leq \ker \mathbf j$ implies $j_1 = \dotsb = j_m$ and we have
\begin{equation*}
 \lim_{n \to \infty} n^{-|\sigma|}\sum_{\substack{1 \leq i_1,\dotsc,i_m \leq n\\ \sigma \leq \ker \mathbf i}} \widetilde{\rho}_{(j_1-1)\cdot n + i_1}(\beta_1)\dotsb \widetilde{\rho}_{(j_m -1)\cdot n + i_m}(\beta_m) = \lim_{n \to \infty} \frac{1}{n}\sum_{i=1}^n \widetilde{\rho}_{(j_1-1)\cdot n + i}(\beta_1\beta_2 \dotsb \beta_m).
\end{equation*}
By Lemma \ref{ergodic}, this converges in $|\;|_2$ to
\begin{equation*}
 E[\widetilde\rho_1(\beta_1\beta_2 \dotsb \beta_m)] = E^{(\sigma)}[\widetilde{\rho_1}(\beta_1) \otimes \dotsb \otimes \widetilde\rho_1(\beta_m)].
\end{equation*}

Now let $\sigma \in NC(m)$ and let $V = \{l+1,\dotsc,l+s\}$ be an interval of $\sigma$,  and let $j$ be the common value of $j_{l+1},\dotsc,j_{l+s}$.  We have
\begin{multline*}
n^{-|\sigma|}\sum_{\substack{1 \leq i_1,\dotsc,i_m \leq n\\ \sigma \leq \ker \mathbf i}} \widetilde{\rho}_{(j_1-1)\cdot n + i_1}(\beta_1)\dotsb \widetilde{\rho}_{(j_m -1)\cdot n + i_m}(\beta_m) \\
= n^{-|\sigma\setminus V|}\negthickspace\sum_{\substack{1 \leq i_1,\dotsc,i_l,\\ i_{l+s+1},\dotsc,i_m \leq n\\ \sigma \setminus V \leq \ker \mathbf i}} \widetilde{\rho}_{(j_1-1)\cdot n + i_1}(\beta_1)\dotsb \Bigl(\frac{1}{n}\sum_{i=1}^n \widetilde{\rho}_{(j-1)\cdot n + i}(\beta_{l+1}\dotsb \beta_{l+s})\Bigr)\dotsb \widetilde{\rho}_{(j_m -1)\cdot n + i_m}(\beta_m) 
\end{multline*}
As above, the interior sum converges to $E[\widetilde \rho_1(\beta_{l+1}\dotsb \beta_{l+s})]$ in $|\;|_2$ as $n \to \infty$.  Now for any $\beta \in C* B$, since the variables $\widetilde{\rho}_i(\beta)$ are identically $*$-distributed with respect to the faithful trace $\tau$, it follows that $\|\widetilde{\rho}_i(\beta)\|$ is independent of $i$.  Therefore there is a constant $D$ such that
\begin{equation*}
 |\widetilde\rho_{i_1}(\beta_1)\dotsb \widetilde\rho_{i_l}(\beta_l) \cdot \xi \cdot \widetilde\rho_{i_{l+s+1}}(\beta_{l+s+1})\dotsb \widetilde \rho_{i_m}(\beta_m)|_2 \leq D|\xi|_2
\end{equation*}
for any $\xi \in L^2(M)$ and $i_1,\dotsc,i_m \in \N$.  It follows that
\begin{multline*}
 \lim_{n \to \infty} n^{-|\sigma\setminus V|}\sum_{\substack{1 \leq i_1,\dotsc,i_l,\\ i_{l+s+1},\dotsc,i_m \leq n\\ \sigma \setminus V \leq \ker \mathbf i}} \widetilde{\rho}_{(j_1-1)\cdot n + i_1}(\beta_1)\dotsb \Bigl(\frac{1}{n}\sum_{i=1}^n \widetilde{\rho}_{j}(\beta_{l+1}\dotsb \beta_{l+s})\Bigr)\dotsb \widetilde{\rho}_{(j_m -1)\cdot n + i_m}(\beta_m) \\
= \lim_{n \to \infty} n^{-|\sigma\setminus V|}\sum_{\substack{1 \leq i_1,\dotsc,i_l,\\ i_{l+s+1},\dotsc,i_m \leq n\\ \sigma \setminus V \leq \ker \mathbf i}} \widetilde{\rho}_{(j_1-1)\cdot n + i_1}(\beta_1)\dotsb E[\widetilde \rho_1(\beta_{l+1}\dotsb \beta_{l+s})]\dotsb \widetilde{\rho}_{(j_m -1)\cdot n + i_m}(\beta_m).
\end{multline*}
By induction, this converges in $|\;|_2$ to
\begin{equation*}
 E^{(\sigma \setminus V)}[\widetilde \rho_1(\beta_1) \otimes \dotsb \otimes \widetilde \rho_1(\beta_l)\cdot E[\widetilde \rho_1(\beta_{l+1}\dotsb \beta_{l+s})] \otimes \dotsb \otimes \widetilde \rho_1(c_m)],
\end{equation*} 
which is precisely $E^{(\sigma)}[\widetilde \rho_1(\beta_1) \otimes \dotsb \otimes \widetilde \rho_1(\beta_m)]$, as desired.

\end{proof}

We can now complete the proof of Theorem \ref{mainthm}.

\begin{proof}[Proof of (ii)$\Rightarrow$(iii)]
Fix $\beta_1,\dotsc,\beta_m \in C * B$ and $1 \leq j_1,\dotsc,j_m \leq k$.  By Proposition \ref{expspread}, for each $n \in \N$ we have
\begin{equation*}
 E[\widetilde{\rho}_{j_1}(\beta_1)\dotsb \widetilde{\rho}_{j_m}(\beta_m)] \otimes 1_{A_i(k,k\cdot n)} = \sum_{1 \leq l_1,\dotsc,l_m \leq kn} E[\widetilde{\rho}_{l_1}(\beta_1)\dotsb \widetilde{\rho}_{l_m}(\beta_m)] \otimes u_{l_1j_1}\dotsb u_{l_mj_m}.
\end{equation*}
Applying $(\mathrm{id} \otimes \psi_{k,n})$, with $\psi_{k,n}$ from Proposition \ref{measure}, to each side of the above equation, we obtain
\begin{multline*}
 E[\widetilde{\rho}_{j_1}(\beta_1)\dotsb \widetilde{\rho}_{j_m}(\beta_m)] \\
= \sum_{1 \leq i_1,\dotsc,i_m \leq n} E[\widetilde{\rho}_{(j_1-1)\cdot n + i_1}(\beta_1)\dotsb \widetilde{\rho}_{(j_m-1)\cdot n + i_m}(\beta_m)] \sum_{\substack{\pi \in NC(m)\\ \pi \leq \ker \mathbf j}} \sum_{\substack{\sigma \in NC(m) \\ \sigma \leq \pi \wedge \ker \mathbf i}} \mu_m(\sigma,\pi)n^{-|\sigma|}\\
= \sum_{\substack{\pi \in NC(m)\\ \pi \leq \ker \mathbf j}} \sum_{\substack{\sigma \in NC(m)\\ \sigma \leq \pi}}\mu_m(\sigma,\pi)E\Bigl[n^{-|\sigma|} \sum_{\substack{1 \leq i_1,\dotsc,i_m \leq n \\ \sigma \leq \ker \mathbf i}} \widetilde{\rho}_{(j_1-1)\cdot n + i_1}(\beta_1)\dotsb \widetilde{\rho}_{(j_m-1)\cdot n + i_m}(\beta_m)\Bigr].
\end{multline*}
Letting $n \to \infty$ and applying Proposition \ref{expform}, we have
\begin{align*}
E[\widetilde{\rho}_{j_1}(\beta_1)\dotsb \widetilde{\rho}_{j_m}(\beta_m)] &= \sum_{\substack{\pi \in NC(m)\\ \pi \leq \ker \mathbf j}} \sum_{\substack{\sigma \in NC(m)\\ \sigma \leq \pi}}\mu_m(\sigma,\pi) E^{(\sigma)}[\widetilde{\rho}_1(\beta_1) \otimes \dotsb \otimes \widetilde{\rho}_1(\beta_m)]\\
&= \sum_{\substack{\pi \in NC(m)\\ \pi \leq \ker \mathbf j}} \kappa_E^{(\pi)}[\widetilde{\rho}_1(\beta_1) \otimes \dotsb \otimes \widetilde{\rho}_1(\beta_m)],
\end{align*}
and the result now follows from Corollary \ref{vancum}.
\end{proof}

\begin{acknowledgement}
I would like to thank my thesis advisor, Dan-Virgil Voiculescu, for his continued guidance and support while completing this project.
\end{acknowledgement}

\def\cprime{$'$}


\begin{thebibliography}{10}

\bibitem{bbc}
{\em T.~Banica, J.~Bichon, and B.~Collins}, Quantum permutation groups: a
  survey, in Noncommutative harmonic analysis with applications to probability,
  vol.~78 of Banach Center Publ., Polish Acad. Sci. Inst. Math., Warsaw, 2007,
  13--34.

\bibitem{bc2}
{\em T.~Banica and B.~Collins}, Integration over quantum permutation groups, J.
  Funct. Anal., {\bf 242} (2007), 641--657.

\bibitem{bcs2}
{\em T.~Banica, S.~Curran, and R.~Speicher}, De {F}inetti theorems for easy
  quantum groups, Ann. Probab., to appear. \href{http://arxiv.org/abs/0907.3314}{{\tt arXiv:0907.3314
  [math.OA]}}, 2009.

\bibitem{cur3}
{\em S.~Curran}, Quantum exchangeable sequences of algebras, Indiana Univ.
  Math. J. {\bf 58} (2009), 1097--1126.

\bibitem{cur4}
{\em S.~Curran}, Quantum rotatability, Trans. Amer. Math. Soc. {\bf 362}
  (2010), 4831--4851.

\bibitem{hs}
{\em E.~Hewitt and L.~J. Savage}, Symmetric measures on {C}artesian products,
  Trans. Amer. Math. Soc. {\bf 80} (1955), 470--501.

\bibitem{kal2}
{\em O.~Kallenberg}, Spreading-invariant sequences and processes on bounded
  index sets, Probab. Theory Related Fields {\bf 118} (2000), 211--250.

\bibitem{kal}
{\em O.~Kallenberg}, {Probabilistic symmetries and invariance principles},
  Probability and its Applications, Springer, New York, 2005.

\bibitem{kos}
{\em C.~K{\"o}stler}, A noncommutative extended de {F}inetti theorem, J. Funct.
  Anal. {\bf 258} (2010), 1073--1120.

\bibitem{ksp}
{\em C.~K{\"o}stler and R.~Speicher}, A noncommutative de {F}inetti theorem:
  invariance under quantum permutations is equivalent to freeness with
  amalgamation, Comm. Math. Phys. {\bf 291} (2009), 473--490.

\bibitem{ns}
{\em A.~Nica and R.~Speicher}, {Lectures on the combinatorics of free
  probability}, vol.~335 of London Mathematical Society Lecture Note Series,
  Cambridge University Press, Cambridge, 2006.

\bibitem{rn}
{\em C.~Ryll-Nardzewski}, On stationary sequences of random variables and the
  de {F}inetti's equivalence, Colloq. Math. {\bf 4} (1957), 149--156.

\bibitem{soltan}
{\em P.~M. So{\l}tan}, Quantum families of maps and quantum semigroups on
  finite quantum spaces, J. Geom. Phys. {\bf 59} (2009), 354--368.

\bibitem{sp1}
{\em R.~Speicher}, Combinatorial theory of the free product with amalgamation
  and operator-valued free probability theory, Mem. Amer. Math. Soc., {\bf 132}
  (1998), x+88.

\bibitem{voi0}
{\em D.~Voiculescu}, Symmetries of some reduced free product
  {$C^\ast$}-algebras, in Operator algebras and their connections with topology
  and ergodic theory, vol.~1132 of Lecture Notes in Math.,
  Springer, Berlin, 1985, 556--588.

\bibitem{vdn}
{\em D.~Voiculescu, K.~Dykema, and A.~Nica}, {Free random variables},
  vol.~1 of CRM Monograph Series, American Mathematical Society, Providence,
  RI, 1992.

\bibitem{wang2}
{\em S.~Wang}, Quantum symmetry groups of finite spaces, Comm. Math. Phys.
  {\bf 195} (1998), 195--211.

\bibitem{wor1}
{\em S.~L. Woronowicz}, Compact matrix pseudogroups, Comm. Math. Phys. {\bf
  111} (1987), 613--665.

\end{thebibliography}
\end{document}